\documentclass[secnum,leqno]{rims-bessatsu}
\usepackage{amssymb, amsmath}
\usepackage{bm,bbm}
\usepackage[all,cmtip]{xy}
\usepackage{graphics}
\usepackage[dvips]{graphicx}
\usepackage{eepic}
\usepackage{epic}
\newtheorem{theorem}{Theorem}[section]

\newtheorem{proposition}[theorem]{Proposition}
\newtheorem{corollary}[theorem]{Corollary}

\theoremstyle{definition}
\newtheorem{definition}[theorem]{Definition}

\newtheorem{remark}[theorem]{Remark}
\newtheorem{question}[theorem]{Question}

\theoremstyle{remark}

\addtocounter{page}{28}   
\newcommand{\bZ}{\mathbbm{Z}}
\newcommand{\bQ}{\mathbbm{Q}}
\newcommand{\bC}{\mathbbm{C}}

\newcommand{\bk}{\mathbbm{k}}\newcommand{\bF}{\mathbbm{F}}

\def\bm#1{\mathbbm{#1}}

\title{Noether's problem and rationality problem for multiplicative invariant fields: a survey}
\TitleHead{Noether's problem and rationality problem for multiplicative invariant fields}          
\author{Akinari \textsc{Hoshi}\footnote{Department of Mathematics, Niigata University, Niigata 950-2181, Japan.\newline e-mail: \texttt{hoshi@math.sc.niigata-u.ac.jp}}}
\AuthorHead{Akinari Hoshi}         
\support{This work was partially supported by JSPS KAKENHI Grant 
Number 25400027.}  
\VolumeNo{77}           
\YearNo{2020}           
\PagesNo{029--053}      
\communication{Received March 27, 2017.
Revised January 31, 2018.}      
\begin{document}
%

\maketitle

\begin{abstract}      
In this paper, 
we give a brief survey of recent developments 
on Noether's problem and rationality problem 
for multiplicative invariant fields 
including 
author's recent papers 
Hoshi \cite{Hos15} 
about Noether's problem over $\bQ$, 
Hoshi, Kang and Kunyavskii \cite{HKK13}, 
Chu, Hoshi, Hu and Kang \cite{CHHK15}, 
Hoshi \cite{Hos16} and 
Hoshi, Kang and Yamasaki \cite{HKY16} 
about Noether's problem over $\bC$, 
and 
Hoshi, Kang and Kitayama \cite{HKK14} and 
Hoshi, Kang and Yamasaki \cite{HKY} 
about rationality problem 
for multiplicative invariant fields. 
\end{abstract}
\tableofcontents      

\section{Introduction}
Let $k$ be a field and 
$G$ be a finite group acting on the rational function field 
$k(x_g \mid g\in G)$ by $k$-automorphisms $h(x_g)=x_{hg}$ for any $g,h\in G$. 
We denote the fixed field ${k(x_g \mid g\in G)^G}$ by $k(G)$. 
Emmy Noether \cite{Noe13, Noe17} asked 
whether $k(G)$ is rational (= purely transcendental) over $k$. 
This is called Noether's problem for $G$ over $k$, 
and is related to the inverse Galois problem,
to the existence of generic $G$-Galois extensions over $k$, and
to the existence of versal $G$-torsors over $k$-rational field extensions
(see Swan \cite{Swa81, Swa83}, Saltman \cite{Sal82}, 
Manin and Tsfasman \cite{MT86}, 
Garibaldi, Merkurjev and Serre \cite[Section 33.1, page 86]{GMS03}, 
Colliot-Th\'el\`ene and Sansuc \cite{CTS07}). 
\begin{theorem}[{Fischer \cite{Fis15}, see also Swan \cite[Theorem 6.1]{Swa83}}]\label{thFis}
Let $G$ be a finite abelian group with exponent $e$. 
Assume that {\rm (i)} either {\rm char} $k=0$ or {\rm char} $k>0$ with 
{\rm char} $k$ ${\not |}$ $e$, and 
{\rm (ii)} $k$ contains a primitive $e$-th root of unity. 
Then $k(G)$ is rational over $k$. 
In particular, $\bm{C}(G)$ is rational over $\bm{C}$.
\end{theorem}
\begin{theorem}[{Kuniyoshi \cite{Kun54, Kun55, Kun56}, see also Gasch\"utz \cite{Gas59}}]\label{thKun}\\
Let $G$ be a $p$-group and $k$ be a field with {\rm char} $k=p>0$. 
Then $k(G)$ is rational over $k$. 
\end{theorem}
%


\begin{definition}\label{defr}
Let $K/k$ and $L/k$ be finitely generated extensions of fields.\\ 
{\rm (1)} $K$ is said to be {\it rational} over $k$ 
(for short, $k$-{\it rational}) if $K$ is purely transcendental over $k$, 
i.e. $K\simeq k(x_1,\ldots,x_n)$ for some algebraically independent elements 
$x_1,\ldots,x_n$ over $k$;\\
{\rm (2)} $K$ is said to be {\it stably $k$-rational} if $K(y_1,\ldots,y_m)$ is 
$k$-rational for some algebraically independent elements 
$y_1,\ldots,y_m$ over $K$;\\ 
{\rm (3)} $K$ and $L$ are said to be {\it stably $k$-isomorphic} 
if $K(y_1,\ldots,y_m)\simeq L(z_1,\ldots,z_n)$ 
for some algebraically independent elements 
$y_1,\ldots,y_m$ over $K$ and $z_1,\ldots,z_n$ over $L$;\\
{\rm (4)} (Saltman, \cite[Definition 3.1]{Sal84b}) 
when $k$ is an infinite field, 
$K$ is said to be 
{\it retract $k$-rational} if there exists a $k$-algebra
$A$ contained in $K$ such that (i) $K$ is the quotient field of
$A$, (ii) there exist a non-zero polynomial $f\in
k[x_1,\ldots,x_n]$ and $k$-algebra homomorphisms $\varphi\colon A\to
k[x_1,\ldots,x_n][1/f]$ and $\psi\colon k[x_1,\ldots,x_n][1/f]\to
A$ satisfying $\psi\circ\varphi =1_A$;\\ 
{\rm (5)} $K$ is said to be {\it $k$-unirational} 
if $k\subset K\subset k(x_1,\ldots,x_n)$ for some integer $n$. 
\end{definition}
We see that 
if $K$ and $L$ are stably $k$-isomorphic and $K$ is retract $k$-rational, 
then $L$ is also retract $k$-rational (see \cite[Proposition 3.6]{Sal84b}), 
and hence it is not difficult to verify the following implications:
\begin{center}
$k$-rational\ \ $\Rightarrow$\ \ 
stably $k$-rational\ \ $\Rightarrow$\ \ 
retract $k$-rational\ \ $\Rightarrow$\ \ 
$k$-unirational. 
\end{center}
Note that $k(G)$ is retract $k$-rational 
if and only if there exists a generic $G$-Galois extension over $k$ 
(see \cite[Theorem 5.3]{Sal82}, \cite[Theorem 3.12]{Sal84b}). 
In particular, if $k$ is a Hilbertian field, e.g. number field, 
and $k(G)$ is retract $k$-rational, then inverse Galois problem 
for $G$ over $k$ has a positive answer, i.e. 
there exists a Galois extension $K/k$ with ${\rm Gal}(K/k)\simeq G$.

\section{Noether's problem over $\bQ$}

Masuda \cite{Mas55, Mas68} gave an idea to use a technique of Galois descent 
to Noether's problem for cyclic groups $C_p$ of order $p$. 
Let $\zeta_p$ be a primitive $p$-th root of unity, 
$L=\bQ(\zeta_p)$ and $\pi={\rm Gal}(L/\bQ)$. 
Then, by Theorem \ref{thFis}, we have $\bQ(C_p)=\bQ(x_1,\ldots,x_p)^{C_p}
=(L(x_1,\ldots,x_p)^{C_p})^{\pi}
=L(y_0,\ldots,y_{p-1})^{\pi}=L(M)^{\pi}(y_0)$ 
where 
$y_0=\sum_{i=1}^p x_i$ is $\pi$-invariant, 
$M$ is free $\bZ[\pi]$-module and $\pi$ acts on $y_1,\ldots,y_{p-1}$ by 
$\sigma(y_i)=\prod_{j=1}^{p-1}y_j^{a_{ij}}$, 
$[a_{ij}]\in GL_n(\bZ)$ for any $\sigma\in\pi$. 
Thus the field $L(M)^\pi$ may be regarded as the function field 
of some algebraic torus of dimension $p-1$ 
(see e.g. \cite[Chapter 3]{Vos98}, \cite[Chapter 1]{HY17}). 

\begin{theorem}[{Masuda \cite{Mas55, Mas68}, see also \cite[Lemma 7.1]{Swa83}}]\label{thMasuda}~\\
{\rm (1)} $M$ is projective $\bZ[\pi]$-module of rank one;\\
{\rm (2)} 
If $M$ is a permutation $\bZ[\pi]$-module, i.e. 
$M$ has a $\bZ$-basis which is permuted by $\pi$, 
then $L(M)^{\pi}$ is $\bQ$-rational. 
In particular, $\bQ(C_p)$ is $\bQ$-rational 
for $p\leq 11$.\footnote{The author \cite[Chapter 5]{Hos05} 
generalized Theorem \ref{thMasuda} (2) to Frobenius groups $F_{pl}$ of order $pl$ 
with $l\mid p-1$ ($p\leq 11$).}
\end{theorem}
Swan \cite{Swa69} gave the first negative solution to Noether's problem 
by investigating a partial converse to Masuda's result. 
\begin{theorem}[{Swan \cite{Swa69}, Voskresenskii \cite{Vos70}}]\label{thSwan}
~\\
{\rm (1)} If $\bQ(C_p)$ is $\bQ$-rational, 
then there exists $\alpha\in \bZ[\zeta_{p-1}]$ such that 
$N_{\bQ(\zeta_{p-1})/\bQ}(\alpha)=\pm p$;\\
{\rm (2)} {\rm (Swan \cite[Theorem 1]{Swa69})} $\bQ(C_{47})$, $\bQ(C_{113})$ and 
$\bQ(C_{233})$ are not $\bQ$-rational;\\
{\rm (3)} {\rm (Voskresenskii \cite[Theorem 2]{Vos70})} $\bQ(C_{47})$, $\bQ(C_{167})$, 
$\bQ(C_{359})$, $\bQ(C_{383})$, $\bQ(C_{479})$, $\bQ(C_{503})$ and 
$\bQ(C_{719})$ are not $\bQ$-rational.
\end{theorem}
\begin{theorem}[{Voskresenskii \cite[Theorem 1]{Vos71}}]\label{thVos}
$\bQ(C_p)$ is $\bQ$-rational if and only if 
there exists $\alpha\in \bZ[\zeta_{p-1}]$ such that 
$N_{\bQ(\zeta_{p-1})/\bQ}(\alpha)=\pm p$.
\end{theorem}
Hence if the cyclotomic field $\bQ(\zeta_{p-1})$ has class number one, 
then $\bQ(C_p)$ is $\bQ$-rational. 
However, it is known that such primes are exactly $p\leq 43$ and $p=61,67,71$ 
(see Masley and Montgomery \cite[Main theorem]{MM76} or 
Washington's book \cite[Chapter 11]{Was97}). 

Endo and Miyata \cite{EM73} refined Masuda-Swan's method 
and gave some further consequences on Noether's problem when $G$ is abelian 
(see also \cite{Vos73}). 
\begin{theorem}[Endo and Miyata {\cite[Theorem 2.3]{EM73}}]\label{thEM23}
Let $G_1$ and $G_2$ be finite groups and $k$ be a field with {\rm char} $k=0$. 
If $k(G_1)$ and $k(G_2)$ are $k$-rational $($resp. stably $k$-rational$)$, 
then 
$k(G_1\times G_2)$ is $k$-rational $($resp. stably $k$-rational$)$.\footnote{Kang and Plans \cite[Theorem 1.3]{KP09} showed that Theorem \ref{thEM23} is also valid for any field $k$.}
\end{theorem}
%
%
\begin{theorem}[Endo and Miyata {\cite[Theorem 3.1]{EM73}}]
Let $p$ be an odd prime and $l$ be a positive integer. 
Let $k$ be a field with {\rm char} $k=0$ and 
$[k(\zeta_{p^l}):k]=p^{m_0}d_0$ with $0\leq m_0\leq l-1$ and 
$d_0\mid p-1$. 
Then the following conditions are equivalent:\\
{\rm (1)} For any faithful $k[C_{p^l}]$-module $V$, 
$k(V)^{C_{p^l}}$ is $k$-rational;\\
{\rm (2)} $k(C_{p^l})$ is $k$-rational;\\
{\rm (3)} There exists $\alpha\in \bZ[\zeta_{p^{m_0}d_0}]$ 
such that 
\[
N_{\bQ(\zeta_{p^{m_0}d_0})/\bQ}(\alpha)=
\begin{cases}
\pm p& m_0>0\\
\pm p^l & m_0=0.
\end{cases}
\]
Further suppose that $m_0>0$. 
Then the above conditions are equivalent to each of 
the following conditions:\\
{\rm (1${}^\prime$)} For any $k[C_{p^l}]$-module $V$, 
$k(V)^{C_{p^l}}$ is $k$-rational;\\
{\rm (2${}^\prime$)} For any $1\leq l^\prime\leq l$, 
$k(C_{p^{l^\prime}})$ is $k$-rational.
\end{theorem}
\begin{theorem}[Endo and Miyata {\cite[Proposition 3.2]{EM73}}]
Let $p$ be an odd prime and $k$ be a field with {\rm char} $k=0$. 
If $k$ contains $\zeta_p+\zeta_p^{-1}$, then $k(C_{p^l})$ is 
$k$-rational for any $l$. 
In particular, $\bQ(C_{3^l})$ is $\bQ$-rational for any $l$. 
\end{theorem}
\begin{theorem}[Endo and Miyata {\cite[Proposition 3.4, Corollary 3.10]{EM73}}]
~\\
{\rm (1)} For primes $p\leq 43$ and $p=61,67,71$, 
$\bQ(C_p)$ is $\bQ$-rational;\\
{\rm (2)} For $p=5,7$, $\bQ(C_{p^2})$ is $\bQ$-rational;\\
{\rm (3)} For $l\geq 3$, $\bQ(C_{2^l})$ is not stably $\bQ$-rational. 
\end{theorem}
\begin{theorem}[Endo and Miyata {\cite[Theorem 4.4]{EM73}}]\label{thEM73m}
Let $G$ be a finite abelian group of odd order and $k$ be a 
field with {\rm char} $k=0$. Then there exists an integer $m>0$ 
such that $k(G^m)$ is $k$-rational. 
\end{theorem}
\begin{theorem}[Endo and Miyata {\cite[Theorem 4.6]{EM73}}]
Let $G$ be a finite abelian group. 
Then $\bQ(G)$ is $\bQ$-rational if and only if 
$\bQ(G)$ is stably $\bQ$-rational. 
\end{theorem}
Ultimately, 
Lenstra \cite{Len74} gave a necessary and sufficient condition 
of Noether's problem for abelian groups. 
\begin{theorem}[Lenstra {\cite[Main Theorem, Remark 5.7]{Len74}}]
Let $k$ be a field and $G$ be a finite abelian group. 
Let $k_{\rm cyc}$ be the maximal cyclotomic extension of $k$ in 
an algebraic closure. 
For $k\subset K\subset k_{\rm cyc}$, we assume that 
$\rho_K={\rm Gal}(K/k)=\langle\tau_k\rangle$ is finite cyclic.
Let $p$ be an odd prime with $p\neq {\rm char}$ $k$ 
and $s\geq 1$ be an integer. 
Let $\mathfrak{a}_K(p^s)$ be a $\bZ[\rho_K]$-ideal defined by 
\begin{align*}
\mathfrak{a}_K(p^s)=
\begin{cases}
\bZ[\rho_K] & {\rm if}\ K\neq k(\zeta_{p^s})\\
(\tau_K-t,p) & {\rm if}\ K=k(\zeta_{p^s})\ {\rm where}\ 
t\in\bZ\ {\rm satisfies}\ \tau_K(\zeta_p)=\zeta_p^t
\end{cases}
\end{align*}
and put $\mathfrak{a}_K(G)=\prod_{p,s}\mathfrak{a}_K(p^s)^{m(G,p,s)}$ 
where $m(G,p,s)=\dim_{\bZ/p\bZ}(p^{s-1}G/p^sG)$. 
Then the following conditions are equivalent:\\
{\rm (1)} $k(G)$ is $k$-rational;\\
{\rm (2)} $k(G)$ is stably $k$-rational;\\
{\rm (3)} for $k\subset K\subset k_{\rm cyc}$, 
the $\bZ[\rho_K]$-ideal $\mathfrak{a}_K(G)$ is principal 
and if {\rm char} $k\neq 2$, then $k(\zeta_{r(G)})/k$ is cyclic extension 
where $r(G)$ is the highest power of $2$ dividing the exponent of $G$. 
\end{theorem}
\begin{theorem}[Lenstra {\cite[Corollary\! 7.2]{Len74}, 
\cite[Proposition\! 2, Corollary\! 3]{Len80}}]\label{thLen74}
Let $n$ be a positive integer. 
Then the following conditions are equivalent:\\
{\rm (1)} $\bQ(C_n)$ is $\bQ$-rational;\\
{\rm (2)} $k(C_n)$ is $k$-rational for any field $k$;\\
{\rm (3)} $\bQ(C_{p^s})$ is $\bQ$-rational 
for any $p^s$ $||$ $n$;\\
{\rm (4)} $8$ ${\not |}$ $n$ and for any $p^s$ $||$ $n$, 
there exists $\alpha\in \bZ[\zeta_{\varphi(p^s)}]$ such that 
$N_{\bQ(\zeta_{\varphi(p^s)})/\bQ}(\alpha)=\pm p$. 
\end{theorem}
\begin{theorem}[{Lenstra \cite[Corollary 7.6]{Len74}, 
\cite[Proposition 6]{Len80}}]
Let $k$\\ be a field which is finitely generated over its prime field. 
Let $P_k$ be the set of primes  $p$ for which 
$k(C_p)$ is $k$-rational. Then $P_k$ has Dirichlet density $0$ 
inside the set of all primes. 
In particular, 
\[
\lim_{x\rightarrow\infty}\frac{\pi^*(x)}{\pi(x)}=0
\]
where $\pi(x)$ is the number of primes $p\leq x$, 
and $\pi^*(x)$ is the number of primes $p\leq x$ 
for which $\bQ(C_p)$ is $\bQ$-rational. 
\end{theorem}
\begin{theorem}[Lenstra {\cite[Proposition 4]{Len80}}]\label{thLen80}
Let $p$ be a prime and $s\geq 2$ be an integer. 
Then $\bQ(C_{p^s})$ is $\bQ$-rational if and only if 
$p^s\in\{2^2, 3^m, 5^2,7^2\mid m\geq 2\}$. 
\end{theorem}
%
By using Theorem \ref{thEM23},
Endo and Miyata \cite[Appendix]{EM73} 
checked whether $\bQ(C_p)$ is $\bQ$-rational for some primes $p<2000$. 
By using PARI/GP \cite{PARI2}, Hoshi \cite{Hos15} 
confirmed that for primes $p<20000$, $\bQ(C_p)$ is not $\bQ$-rational 
except for $17$ rational cases with $p\leq 43$ and $p=61, 67, 71$ 
and undetermined $46$ cases. 
Eventually, Plans \cite{Pla17} determined the complete set of primes 
for which $\bQ(C_p)$ is $\bQ$-rational: 

\begin{theorem}[{Plans \cite[Theorem 1.1]{Pla17}}]\label{thP}
Let $p$ be a prime. Then 
$\bQ(C_p)$ is $\bQ$-rational if and only if 
$p\leq 43$, $p=61, 67$ or $71$. 
\end{theorem}

Combining  Theorem \ref{thLen74}, Theorem \ref{thLen80} 
and Theorem \ref{thP}, we have:  

\begin{corollary}[{Plans \cite[Corollary 1.2]{Pla17}}]
Let $n$ be a positive integer. 
Then $\bQ(C_n)$ is $\bQ$-rational if and only if 
$n$ divides 
\[
2^2\cdot 3^m\cdot 5^2\cdot 7^2\cdot 11\cdot 13\cdot 17\cdot 19\cdot 
23\cdot 29\cdot 31\cdot 37\cdot 41\cdot 43\cdot 61\cdot 67\cdot 71
\]
for some integer $m\geq 0$.
\end{corollary}
On the other hand, just a handful of results about 
Noether's problem are obtained when the groups are non-abelian.  
\begin{theorem}[{Maeda \cite[Theorem, page 418]{Mae89}}]
Let $k$ be a field and $A_5$ be the alternating group of degree $5$. 
Then $k(A_5)$ is $k$-rational. 
\end{theorem}
\begin{theorem}[{Rikuna \cite{Rik}, Plans \cite{Pla07}, see also \cite[Example 13.7]{HKY11}}]\\
Let $k$ be a field with ${\rm char}$ $k\neq 2$. Then 
$k(SL_2(\bF_3))$ and $k(GL_2(\bF_3))$ are $k$-rational. 
\end{theorem}
\begin{theorem}[{Serre \cite[Chapter IX]{GMS03}, see also Kang \cite{Kan05}}]
Let $G$ be a finite group with a $2$-Sylow subgroup which is cyclic of 
order $\geq 8$ or the generalized quaternion $Q_{16}$ of order $16$. 
Then $\bQ(G)$ is not stably $\bQ$-rational. 
\end{theorem}
\begin{theorem}[{Plans \cite[Theorem 2]{Pla09}}]
Let $A_n$ be the alternating group of degree $n$. 
If $n\geq 3$ is odd integer, then 
$\bQ(A_n)$ is rational over $\bQ(A_{n-1})$. 
In particular, if $\bQ(A_{n-1})$ is $\bQ$-rational, then so is $\bQ(A_n)$. 
\end{theorem}
However, it is an open problem whether $k(A_n)$ is $k$-rational 
for $n\geq 6$.

\section{Noether's problem over $\bC$ and unramified Brauer groups}
We consider Noether's problem for $G$ over $\bC$, i.e. 
the rationality problem for $\bC(G)$ over $\bC$. 
Let $G$ be a $p$-group. 
Then, by Theorem \ref{thFis} and Theorem \ref{thKun}, 
we may focus on the case where $G$ is a non-abelian $p$-group 
and $k$ is a field with char $k\neq p$. 
For $p$-groups of small order, the following results are known.
\begin{theorem}[Chu and Kang \cite{CK01}] \label{thCK01}
Let $p$ be any prime and $G$ be a $p$-group of order $\le p^4$
and of exponent $e$. If $k$ is a field containing a primitive $e$-th root
of unity, then $k(G)$ is $k$-rational. 
In particular, $\bC(G)$ is $\bC$-rational. 
\end{theorem}
\begin{theorem}[Chu, Hu, Kang and Prokhorov \cite{CHKP08}]\label{thCHKP08}
Let $G$ be a group of order $32$ and of exponent $e$.  
If $k$ is a field containing a primitive $e$-th root of unity, 
then $k(G)$ is $k$-rational. 
In particular, $\bC(G)$ is $\bC$-rational. 
\end{theorem}
%
%
Saltman introduced a notion of retract $k$-rationality 
(see Definition \ref{defr}) and the unramified Brauer group: 
\begin{definition}[{Saltman \cite[Definition 3.1]{Sal84a}, \cite[page 56]{Sal85}}] \label{defSal}
Let $K/k$ be an extension of fields. 
The {\it unramified Brauer group} ${\rm Br}_{\rm nr}(K/k)$ 
of $K$ over $k$ is defined to be 
\[
{\rm Br}_{\rm nr}(K/k)=\bigcap_R {\rm Image} \{ {\rm Br}(R)\to{\rm Br}(K)\}
\]
where ${\rm Br}(R)\to {\rm Br}(K)$ is the natural map of
Brauer groups and $R$ runs over all the discrete valuation rings $R$ 
such that $k\subset R\subset K$ and $K$ is the quotient field of $R$. 
We 
write just 
${\rm Br}_{\rm nr}(K)$ when the base field $k$ is clear from the context.
\end{definition}
\begin{proposition}[{Saltman \cite{Sal84a}, \cite[Proposition 1.8]{Sal85}, \cite{Sal87}}] \label{propSal}
If $K$ is\\ retract $k$-rational, then
${\rm Br}(k)\overset{\sim}{\longrightarrow}{\rm Br}_{\rm nr}(K)$. 
In particular, if $k$ is an algebraically closed field and $K$ is
retract $k$-rational, then ${\rm Br}_{\rm nr}(K)=0$.
\end{proposition}
\begin{theorem}[{Bogomolov \cite[Theorem 3.1]{Bog88}, Saltman \cite[Theorem 12]{Sal90}}]\label{thBog}
Let $G$ be a finite group and $k$ be an algebraically closed field with 
{\rm char} $k=0$ or {\rm char} $k=p$ 
${\not |}$ $|G|$. 
Then ${\rm Br}_{\rm nr}(k(G)/k)$ is isomorphic to the group $B_0(G)$ defined by
\[
B_0(G)=\bigcap_A {\rm Ker}\{{\rm res}: H^2(G,\bQ/\bZ)\to H^2(A,\bQ/\bZ)\}
\]
where $A$ runs over all the bicyclic subgroups of $G$ $($a group $A$
is called bicyclic if $A$ is either a cyclic group or a direct
product of two cyclic groups$)$.
\end{theorem}
\begin{remark}
For a smooth projective variety $X$ over $\bC$ with function 
field $K$, ${\rm Br}_{\rm nr}(K/\bC)$ is isomorphic to 
the birational invariant 
$H^3(X,\bZ)_{\rm tors}$ which was used by Artin and Mumford \cite{AM72} 
to provide some elementary examples of $k$-unirational varieties 
which are not $k$-rational (see also \cite[Theorem 1.1 and Corollary]{Bog88}). 
\end{remark}
Note that $B_0(G)$ is a subgroup of $H^2(G,\bQ/\bZ)$ 
which is isomorphic to the Schur multiplier 
$H_2(G,\bZ)$ of $G$ (see Karpilovsky \cite{Kar87}). 
We call $B_0(G)$ the {\it Bogomolov multiplier of $G$} 
(cf. Kunyavskii \cite{Kun10}). 
Because of Theorem \ref{thBog}, we will not 
distinguish $B_0(G)$ and ${\rm Br}_{\rm nr}(k(G)/k)$ when $k$ is 
an algebraically closed field, and {\rm char} $k=0$ or 
{\rm char} $k=p$ 
${\not |}$ $|G|$. 
Using $B_0(G)$, Saltman and Bogomolov gave 
counter-examples to Noether's problem for non-abelian
$p$-groups over algebraically closed field.

\begin{theorem}[{Saltman \cite{Sal84a}, Bogomolov \cite{Bog88}}] \label{thSB}
Let $p$ be any prime and $k$ be any algebraically closed field
with ${\rm char}$ $k\ne p$.\\
{\rm (1) (Saltman \cite[Theorem 3.6]{Sal84a})} 
There exists a meta-abelian group $G$ of order $p^9$
such that $B_0(G)\ne 0$. 
In particular, $k(G)$ is not $($retract, stably$)$ $k$-rational;\\
{\rm (2) (Bogomolov \cite[Lemma 5.6]{Bog88})} 
There exists a group $G$ of order $p^6$ such that $B_0(G)\ne 0$. 
In particular, $k(G)$ is not $($retract, stably$)$ $k$-rational.
\end{theorem}
Colliot-Th\'el\`ene and Ojanguren \cite{CTO89} 
generalized the notion of the unramified Brauer group 
${\rm Br}_{\rm nr}(K/k)$ to the unramified cohomology 
$H_{\rm nr}^i(K/k,\mu_n^{\otimes j})$ of degree $i\geq 1$, 
that is $F_n^{i,j}(K/k)$ in \cite[Definition 1.1]{CTO89}.
\begin{definition}[{Colliot-Th\'el\`ene and Ojanguren \cite{CTO89}, \cite[Sections 2--4]{CT95}}]\label{d2.6}
Let $n$ be a positive integer and $k$ be a field with {\rm char} $k=0$ or {\rm char} $k=p$ with $p {\not |}$ $n$.
Let $K/k$ be a function field, that is finitely generated field extension as a field over $k$. For any positive integer $i\ge 2 $, any integer $j$, the {\it unramified cohomology group} $H^i_{\rm nr}(K/k,\mu_n^{\otimes j})$
of $K$ over $k$ of degree $i$ is defined to be 
\[
H^i_{\rm nr}(K/k,\mu_n^{\otimes j}):=\bigcap_R {\rm Ker}
\{r_R: H^i(K,\mu_n^{\otimes j})\to H^{i-1}(\bk_R,\mu_n^{\otimes (j-1)})\}
\]
where $R$ runs over all the discrete valuation rings $R$ of rank one
such that $k\subset R\subset K$ and $K$ is the quotient field of $R$, 
$\bk_R$ is the residue field of $R$ 
and $r_R$ 
is the residue map of $K$ at $R$. 
\end{definition}

By \cite[Theorem 4.1.1, page 30]{CT95}, if it is assumed furthermore that $K$ is the function field of a complete smooth variety over $k$, the unramified cohomology group $H^i_{\rm nr}(K/k,\mu_n^{\otimes j})$ may be defined as well by
\[
H^i_{\rm nr}(K/k,\mu_n^{\otimes j})=\bigcap_R {\rm Image}
\{ H^i_{\rm \acute{e}t}(R,\mu_n^{\otimes j})\to H^i_{\rm \acute{e}t}(K,\mu_n^{\otimes j})\}
\]
where $R$ runs over all the discrete valuation rings $R$ of rank one
such that $k\subset R\subset K$ and $K$ is the quotient field of $R$.

Note that the unramified cohomology groups of degree two are isomorphic to 
the $n$-torsion part of the unramified Brauer group: 
${}_n{\rm Br}_{\rm nr}(K/k)\simeq H_{\rm nr}^2(K/k,\mu_n)$.
\begin{theorem} 
Let $n$ be a positive integer and 
$k$ be an algebraically closed field 
with {\rm char} $k=0$ or {\rm char} 
$k=p {\not |}$ $n$.\\
{\rm (1) (Colliot-Th\'el\`ene and Ojanguren \cite[Proposition 1.2]{CTO89})} 
If $K$ and $L$ are stably $k$-isomorphic, then 
$H_{\rm nr}^i(K/k,\mu_n^{\otimes j}) \overset{\sim}{\longrightarrow} H_{\rm nr}^i(L/k,\mu_n^{\otimes j})$.
In particular, 
$K$ is stably $k$-rational, then $H_{\rm nr}^i(K/k,\mu_n^{\otimes j})=0$;\\
{\rm (2) (\cite[Proposition 2.15]{Mer08}, 
see also \cite[Remarque 1.2.2]{CTO89}, \cite[Sections 2--4]{CT95}, 
\cite[Example 5.9]{GS10})} 
If $K$ is retract $k$-rational, then $H_{\rm nr}^i(K/k,\mu_n^{\otimes j})=0$.
\end{theorem}
Colliot-Th\'el\`ene and Ojanguren 
\cite[Section 3]{CTO89} 
produced the first example of not stably $\bC$-rational 
but $\bC$-unirational field $K$ 
with $H_{\rm nr}^3(K,\mu_2^{\otimes 3})\neq 0$, 
where $K$ is the function field of 
a quadric of the type $\langle\!\langle f_1,f_2\rangle\!\rangle
=\langle g_1g_2\rangle$
over the rational function field $\bC(x,y,z)$ 
with three variables $x,y,z$ 
for a $2$-fold Pfister form $\langle\!\langle f_1,f_2\rangle\!\rangle$, 
as a generalization of Artin and Mumford \cite{AM72}. 
Peyre \cite[Corollary 3]{Pey93} gave a sufficient condition for 
$H_{\rm nr}^i(K/k,\mu_p^{\otimes i})\neq 0$ and 
produced an example of the function field $K$ with 
$H_{\rm nr}^3(K/k,\mu_p^{\otimes 3})\neq 0$ and ${\rm Br}_{\rm nr}(K/k)=0$ 
using a result of Suslin \cite{Sus91} 
where $K$ is the function field of 
a product of some norm varieties 
associated to cyclic central simple algebras of degree $p$ 
(see \cite[Proposition 7]{Pey93}). 
Using a result of Jacob and Rost \cite{JR89}, 
Peyre \cite[Proposition 9]{Pey93} 
also gave an example of $H_{\rm nr}^4(K/k,\mu_2^{\otimes 4})\neq 0$ 
and ${\rm Br}_{\rm nr}(K/k)=0$ 
where $K$ is the function field of a product of quadrics associated to 
a $4$-fold Pfister form $\langle\!\langle a_1,a_2,a_3,a_4\rangle\!\rangle$ 
(see also \cite[Section 4.2]{CT95}).

In case ${\rm char} \, k =0$,
take the direct limit with respect to $n$:
\[
H^i(K/k,\bQ/\bZ(j))=
\lim_{\overset{\longrightarrow}{n}}H^i(K/k,\mu_n^{\otimes j})
\]
and we may define the unramified cohomology group
\[
H_{\rm nr}^i(K/k,\bQ/\bZ(j))
=\bigcap_R {\rm Ker}
\{r_R: H^i(K/k,\bQ/\bZ(j))\to H^{i-1}(\bk_R,\bQ/\bZ(j-1))\}.
\]

We write simply $H^i_{\rm nr}(K,\mu_n^{\otimes j})$ and $H_{\rm nr}^i(K,\bQ/\bZ(j))$ when the base field $k$ is understood.
When $k$ is an algebraically closed field with {\rm char} $k=0$, we will write $H_{\rm nr}^i(K/k,\bQ/\bZ)$ for $H_{\rm nr}^i(K/k,\bQ/\bZ(j))$. 
Then we have ${\rm Br}_{\rm nr}(K/k)\simeq H_{\rm nr}^2(K/k,\bQ/\bZ)$. 

Peyre \cite{Pey08} constructed an example of a field $K$, as $K=\bC(G)$, 
whose unramified Brauer group vanishes, but 
unramified cohomology of degree three does not vanish:
\begin{theorem}[{Peyre \cite[Theorem 3]{Pey08}}] \label{thPeyre}
Let $p$ be any odd prime. 
Then there exists a $p$-group $G$ of order $p^{12}$ such that 
$B_0(G)=0$ and 
$H_{\rm nr}^3(\bC(G),\bQ/\bZ)\ne 0$.
In particular, $\bC(G)$ is not $($retract, stably$)$ $\bC$-rational.
\end{theorem}
The idea of Peyre's proof 
is to find a subgroup $K_{\max}^3/K^3$
of $H_{\rm nr}^3(\mathbbm{C}(G),\mathbbm{Q}/\mathbbm{Z})$ and to show that
$K_{\max}^3/K^3\ne 0$ (see \cite[page 210]{Pey08}). 

Asok \cite{Aso13} generalized Peyre's argument \cite{Pey93} and 
established the following theorem 
for a smooth proper model $X$ (resp. a smooth projective model $Y$) 
of the function field of a product of 
quadrics of the type $\langle\!\langle s_1,\ldots,s_{n-1}\rangle\!\rangle
=\langle s_n\rangle$ 
(resp. Rost varieties) 
over some rational function field over $\bC$ with many variables. 
\begin{theorem}[{Asok \cite{Aso13}, see \cite[Theorem 3]{AM11} for retract rationality}]\label{thAsok}\\
{\rm (1)} {\rm (\cite[Theorem 1]{Aso13})} 
For any $n>0$, 
there exists a smooth projective complex variety $X$ that is 
$\bC$-unirational, for which 
$H_{\rm nr}^i(\bC(X),\mu_2^{\otimes i})=0$ for each $i<n$, yet
$H_{\rm nr}^n(\bC(X),\mu_2^{\otimes n})\neq 0$, and so 
$X$ is not $\mathbb{A}^1$-connected, 
nor $($retract, stably$)$ $\bC$-rational;\\
{\rm (2)} {\rm (\cite[Theorem 3]{Aso13})} 
For any prime $l$ and any $n\geq 2$, 
there exists a smooth projective rationally connected complex 
variety $Y$ such that 
$H_{\rm nr}^n(\bC(Y),\mu_l^{\otimes n})\neq 0$. 
In particular, $Y$ is not $\mathbb{A}^1$-connected, 
nor $($retract, stably$)$ $\bC$-rational. 
\end{theorem}
Namely, the triviality of the unramified Brauer group or the unramified
cohomology of higher degree is just a necessary condition of
$\bC$-rationality of fields. 
It is unknown whether the vanishing of
all the unramified cohomologies is a sufficient condition for
$\bC$-rationality. 
It is interesting to consider an analog of Theorem \ref{thAsok} 
for quotient varieties $V/G$, e.g. the case of Noether's problem 
$\bC(V_{\rm reg}/G)=\bC(G)$. 

Colliot-Th\'el\`ene and Voisin \cite{CTV12} established: 
\begin{theorem}[{Colliot-Th\'el\`ene and Voisin {\cite{CTV12}, \cite[Theorem 6.18]{Voi14}}}]\label{thCTV}\\
For any smooth projective complex variety $X$, there is an exact sequence 
\[
0\rightarrow H_{\rm nr}^3(X,\bZ)\otimes\bQ/\bZ\rightarrow 
H^3_{\rm nr}(X,\bQ/\bZ)\rightarrow {\rm Tors}(Z^4(X))\rightarrow 0
\]
where 
\[
Z^4(X)={\rm Hdg}^4(X,\bZ)/{\rm Hdg}^4(X,\bZ)_{\rm alg}
\]
and the lower index ``alg'' means that we consider the group of 
integral Hodge classes which are algebraic. 
In particular, if $X$ is rationally connected, then we have 
\[
H_{\rm nr}^3(X,\bQ/\bZ)\simeq Z^4(X). 
\]
\end{theorem}
Using Peyre's method \cite{Pey08}, we obtain the following theorem 
which is an improvement of Theorem \ref{thPeyre} and 
gives an explicit counter-example to integral Hodge conjecture 
with the aid of Theorem \ref{thCTV}. 
\begin{theorem}[{Hoshi, Kang and Yamasaki \cite[Theorem 1.4]{HKY16}}] \label{t1.4}
Let $p$ be any odd prime. 
Then there exists a $p$-group $G$ of order
$p^9$ such that $B_0(G)=0$ and 
$H_{\rm nr}^3(\mathbbm{C}(G),\mathbbm{Q}/\mathbbm{Z})\ne 0$. 
In particular, $\bC(G)$ is not $($retract, stably$)$ $\bC$-rational.
\end{theorem}
%

{\bf The case where $G$ is a group of order $p^5$ $(p\geq 3)$.}

{}From Theorem \ref{thSB} (2), 
Bogomolov \cite[Remark 1]{Bog88} raised a question 
to classify the groups of order $p^6$ with $B_0(G)\neq 0$. 
He also claimed that if $G$ is a $p$-group of order $\leq p^5$, 
then $B_0(G)=0$ (\cite[Lemma 5.6]{Bog88}). 
However, this claim was disproved by Moravec:
\begin{theorem}[{Moravec \cite[Section 8]{Mor12}}] \label{thMo}
Let $G$ be a group of order $243$. 
Then $B_0(G)\ne 0$ if and only if 
$G=G(3^5,i)$ with $28\le i\le 30$, where $G(3^5,i)$ is the $i$-th
group of order $243$ in the GAP database {\rm \cite{GAP}}. 
Moreover, if $B_0(G)\neq 0$, then $B_0(G)\simeq \bZ/3\bZ$.
\end{theorem}
%
Moravec \cite{Mor12} gave a formula for $B_0(G)$ 
by using a nonabelian exterior square $G \wedge G$ of $G$ 
and an implemented algorithm {\bf b0g.g} in computer algebra system 
GAP \cite{GAP}, which is available from his website 
\verb"www.fmf.uni-lj.si/~moravec/Papers/b0g.g".
The number of all solvable groups $G$ of order $\leq 729$ apart 
from the orders $512$, $576$ and $640$ with $B_0(G)\neq 0$ 
was given as in \cite[Table 1]{Mor12}. 

Hoshi, Kang and Kunyavskii \cite{HKK13} determined $p$-groups $G$ 
of order $p^5$ with $B_0(G)\ne 0$ for any $p\geq 3$. 
It turns out that they belong to the same isoclinism family.
\begin{definition}[{Hall \cite[page 133]{Hal40}}] \label{defHall}
Let $G$ be a finite group. 
Let $Z(G)$ be the center of $G$ and 
$[G,G]$ be the commutator subgroup of $G$. 
Two $p$-groups $G_1$ and $G_2$ are called {\it isoclinic} if there exist
group isomorphisms $\theta\colon G_1/Z(G_1) \to G_2/Z(G_2)$ and
$\phi\colon [G_1,G_1]\to [G_2,G_2]$ such that $\phi([g,h])$
$=[g',h']$ for any $g,h\in G_1$ with $g'\in \theta(gZ(G_1))$, $h'\in
\theta(hZ(G_1))$: 
\[\xymatrix{
G_1/Z_1\times G_1/Z_1 \ar[d]_{[\cdot,\cdot]} \ar[r]^{(\theta,\theta)} 
\ar@{}[dr]| \circlearrowleft 
& G_2/Z_2\times G_2/Z_2 \ar[d]_{[\cdot,\cdot]} \\
[G_1, G_1] \ar[r]^\phi & [G_2, G_2].\\
}\]
For a prime $p$ and an integer $n$, 
we denote by $G_n(p)$ the set of all non-isomorphic groups of order $p^n$. 
In $G_n(p)$, consider an equivalence relation: two groups $G_1$ and $G_2$ are
equivalent if and only if they are isoclinic. 
Each equivalence class of $G_n(p)$ is called an {\it isoclinism family}, 
and the $j$-th isoclinism family is denoted by $\Phi_j$.
\end{definition}
For $p\geq 5$ (resp. $p=3$), there exist $2p+61+\gcd\{4,p-1\}+2\gcd\{3,p-1\}$ 
(resp. $67$) groups $G$ of order $p^5$ which are classified into ten 
isoclinism families $\Phi_1,\ldots,\Phi_{10}$ (see \cite[Section 4]{Jam80}). 
The main theorem of \cite{HKK13} can be stated as follows:
\begin{theorem}[{Hoshi, Kang and Kunyavskii \cite[Theorem 1.12]{HKK13}}] \label{thKKu13}
Let $p$ be any odd prime and $G$ be a group of order $p^5$. Then
$B_0(G)\ne 0$ if and only if $G$ belongs to
the isoclinism family $\Phi_{10}$. 
Moreover, if $B_0(G)\neq 0$, then $B_0(G)\simeq \bZ/p\bZ$.
\end{theorem}
For the last statement, see \cite[Remark, page 424]{Kan14}.
The proof of Theorem \ref{thKKu13} was given by purely algebraic way. 
There exist exactly $3$ groups which belong to $\Phi_{10}$ if $p=3$, i.e. 
$G=G(243,i)$ with $28\leq i\leq 30$. 
This agrees with Moravec's computational result (Theorem \ref{thMo}).
For $p\ge5$, there exist exactly 
$1+\gcd\{4,p-1\}+\gcd \{3,p-1\}$ groups which belong to $\Phi_{10}$ 
(see \cite[page 621]{Jam80}). 

The following result for the $k$-rationality of $k(G)$ 
supplements Theorem \ref{thMo} although it is unknown whether $k(G)$ 
is $k$-rational for groups $G$ which belong to $\Phi_7$: 
\begin{theorem}[{Chu, Hoshi, Hu and Kang \cite[Theorem 1.13]{CHHK15}}]
Let $G$ be a group of order $243$ with exponent $e$. 
If $B_0(G)=0$ and $k$ be a field containing a primitive $e$-th root of unity, 
then $k(G)$ is $k$-rational 
except possibly for the five groups $G$ which belong to $\Phi_7$, 
i.e. $G=G(243,i)$ with $56 \le i \le 60$.
\end{theorem}
%
In \cite{HKK13} and \cite{CHHK15}, 
not only the evaluation of the Bogomolov multiplier 
$B_0(G)$ and the $k$-rationality of $k(G)$ 
but also the $k$-isomorphisms between 
$k(G_1)$ and $k(G_2)$ for some groups $G_1$ and $G_2$ 
belonging to the same isoclinism family were given. 

Bogomolov and B\"ohning \cite{BB13} 
gave an answer to the question raised as 
\cite[Question 1.11]{HKK13} in the affirmative as follows.  
\begin{theorem}[{Bogomolov and B\"ohning \cite[Theorem 6]{BB13}}]\label{thBB}
If $G_1$ and $G_2$ are isoclinic, 
then $\bC(G_1)$ and $\bC(G_2)$ are stably $\bC$-isomorphic. 
In particular, $H_{\rm nr}^i(\bC(G_1),\mu_n^{\otimes j})$ 
$\overset{\sim}{\longrightarrow}$ $H_{\rm nr}^i(\bC(G_2),\mu_n^{\otimes j})$. 
\end{theorem}
A partial result of Theorem \ref{thBB} was already given by Moravec. 
Indeed, 
Moravec \cite[Theorem 1.2]{Mor14} proved that 
if $G_1$ and $G_2$ are isoclinic, then $B_0(G_1)\simeq B_0(G_2)$.\\

{\bf
The case where $G$ is a group of order $64$.}

The classification of the groups $G$ of order $64=2^6$ 
with $B_0(G)\neq 0$ 
was obtained by Chu, Hu, Kang and Kunyavskii \cite{CHKK10}.
Moreover, they investigated Noether's problem for groups $G$ 
with $B_0(G)=0$. 
There exist $267$ groups $G$ of order $64$ which are classified into $27$ 
isoclinism families $\Phi_1,\ldots,\Phi_{27}$ 
by Hall and Senior \cite{HS64} (see also \cite[Table I]{JNO90}). 
The main result of \cite{CHKK10} 
can be stated in terms of the isoclinism families as follows.

\begin{theorem}[{Chu, Hu, Kang and Kunyavskii \cite{CHKK10}}]\label{thCHKK10}
Let $G=G(2^6,i)$, $1\leq i\leq 267$, be the $i$-th group of order $64$ 
in the GAP database {\rm \cite{GAP}}.\\
{\rm (1) (\cite[Theorem 1.8]{CHKK10})} 
$B_0(G)\ne 0$ if and only if $G$ belongs to the isoclinism family 
$\Phi_{16}$, i.e. $G=G(2^6,i)$ with $149\le i\le 151$,
$170\le i\le 172$, $177\le i\le 178$ or $i=182$. 
Moreover, if $B_0(G)\neq 0$, then $B_0(G)\simeq \bZ/2\bZ$ 
$($see {\rm \cite[Remark, page 424]{Kan14}} for this statement$)$;\\
{\rm (2) (\cite[Theorem 1.10]{CHKK10})} 
If $B_0(G)= 0$ and $k$ is an quadratically closed field, 
then $k(G)$ is $k$-rational except
possibly for five groups which belong to $\Phi_{13}$, i.e. 
$G=G(2^6,i)$ with $241\le i\le 245$.
\end{theorem}
For groups $G$ which belong to $\Phi_{13}$, 
$k$-rationality of $k(G)$ is unknown. 
The following two propositions supplement the cases $\Phi_{13}$ 
and $\Phi_{16}$ of Theorem \ref{thCHKK10}. 
For the proof, 
the case of $G=G(2^6,149)$ is given in \cite[Proof of Theorem 6.3]{HKK14}, 
see also \cite[Example 5.11, page 2355]{CHKK10} 
and the proof for other cases can be obtained by the similar manner. 
\begin{definition}
\label{defL01}
Let $k$ be a field with {\rm char} $k\neq 2$ and 
$k(X_1,X_2,X_3,X_4,X_5,X_6)$ be the rational function field over $k$ 
with variables $X_1,X_2,X_3,X_4,X_5,X_6$.\\
{\rm (1)} {\it The field $L_k^{(0)}$} is defined to be $k(X_1,X_2,X_3,X_4,X_5,X_6)^H$ 
where $H=\langle \sigma_1, \sigma_2\rangle\simeq \bZ/2\bZ\oplus\bZ/2\bZ$ 
acts on 
$k(X_1,X_2,X_3,X_4,X_5,X_6)$ by $k$-automorphisms  
\begin{align*}
\sigma_1 &: X_1\mapsto X_3,\ X_2\mapsto \frac{1}{X_1X_2X_3},\ X_3\mapsto X_1,\ 
X_4\mapsto X_6,\ X_5\mapsto \frac{1}{X_4X_5X_6},\ X_6\mapsto X_4,\\
\sigma_2 &: X_1\mapsto X_2,\ X_2\mapsto X_1,\ X_3\mapsto \frac{1}{X_1X_2X_3},\ 
X_4\mapsto X_5,\ X_5\mapsto X_4,\ X_6\mapsto \frac{1}{X_4X_5X_6}. 
\end{align*}
{\rm (2)} {\it The field $L_k^{(1)}$} 
is defined to be $k(X_1,X_2,X_3,X_4)^{\langle\tau\rangle}$ 
where $\langle\tau\rangle \simeq C_2$ acts on 
$k(X_1,X_2,$ $X_3,X_4)$ by $k$-automorphisms 
\begin{align*}
\tau: X_1\mapsto -X_1,\ 
X_2\mapsto \frac{X_4}{X_2},\ 
X_3\mapsto \frac{(X_4-1)(X_4-X_1^2)}{X_3},\ 
X_4\mapsto X_4.
\end{align*}
\end{definition}
\begin{proposition}[{\cite[Proposition 6.3]{CHKK10}, see also \cite[Proposition 12.5]{HY17}}]\label{prop1}
Let $G$ be a group of order $64$ which belongs to $\Phi_{13}$, i.e. 
$G=G(2^6,i)$ with $241\leq i\leq 245$. 
There exists a $\bC$-injective homomorphism $\varphi : L_\bC^{(0)}\rightarrow \bC(G)$ 
such that $\bC(G)$ is rational over $\varphi(L_\bC^{(0)})$. 
In particular, $\bC(G)$ and $L_\bC^{(0)}$ are stably $\bC$-isomorphic 
and $B_0(G)\simeq {\rm Br}_{\rm nr}(L_\bC^{(0)})=0$. 
\end{proposition}
\begin{proposition}[{\cite[Example 5.11]{CHKK10}, \cite[Proof of Theorem 6.3]{HKK14}}]\label{prop2}\\
Let $G$ be a group of order $64$ which belongs to $\Phi_{16}$, 
i.e. $G=G(2^6,i)$ with $149\le i\le 151$,
$170\le i\le 172$, $177\le i\le 178$ or $i=182$. 
There exists a $\bC$-injective homomorphism 
$\varphi : L_\bC^{(1)}\rightarrow \bC(G)$ such that $\bC(G)$ is rational 
over $\varphi(L_\bC^{(1)})$. 
In particular, $\bC(G)$ and $L_\bC^{(1)}$ are stably $\bC$-isomorphic, 
$B_0(G)\simeq {\rm Br}_{\rm nr}(L_\bC^{(1)})\simeq \bZ/2\bZ$ and hence $\bC(G)$ and $L_\bC^{(1)}$ 
are not $($retract, stably$)$ $\bC$-rational. 
\end{proposition}
\begin{question}[{\cite[Section 6]{CHKK10}, \cite[Section 12]{HY17}}]
Is $L_k^{(0)}$ $k$-rational?
\end{question}

\smallskip
{\bf 
The case where $G$ is a group of order $128$.}

There exist $2328$ groups 
of order $128$ which are classified into $115$ 
isoclinism families $\Phi_1,\ldots,\Phi_{115}$ 
(\cite[Tables I, II, III]{JNO90}). 
\begin{theorem}[{Moravec \cite[Section 8, Table 1]{Mor12}}]
Let $G$ be a group of order $128$. Then 
$B_0(G)\neq 0$ if and only if 
$G$ belongs to the isoclinism family 
$\Phi_{16}$, $\Phi_{30}$, 
$\Phi_{31}$, $\Phi_{37}$, $\Phi_{39}$, $\Phi_{43}$, 
$\Phi_{58}$, $\Phi_{60}$, $\Phi_{80}$, $\Phi_{106}$ or $\Phi_{114}$. 
Moreover, 
we have
\begin{align*}
B_0(G)\simeq\begin{cases}
\bZ/2\bZ\!\!\!\!&{\rm if}\ G\ {\rm belongs\ to\ } \Phi_{16}, \Phi_{31}, \Phi_{37}, 
\Phi_{39}, \Phi_{43}, \Phi_{58}, \Phi_{60}, \Phi_{80},\Phi_{106}\ {\rm or}\ 
\Phi_{114},\\
(\bZ/2\bZ)^{\oplus 2}\!\!\!\!&{\rm if}\ G\ {\rm belongs\ to\ } \Phi_{30}.
\end{cases}
\end{align*}
In particular, $\bC(G)$ is not $($retract, stably$)$ $\bC$-rational.
\end{theorem}
It turns out that 
there exist $220$ groups $G$ of order $128$ 
with $B_0(G)\neq 0$: 

\begin{table}[h]
\begin{tabular}{|c|cccccccccc|c|}\hline
Family & $\Phi_{16}$\! &\! $\Phi_{31}$\! &\! $\Phi_{37}$\! &\! $\Phi_{39}$\! &\! $\Phi_{43}$\! &\!\! $\Phi_{58}$\!\! & \! $\Phi_{60}$\! &\! $\Phi_{80}$\!&\!  
$\Phi_{106}$\!&\!$\Phi_{114}$\!&\!$\Phi_{30}$\!\\\hline
$\exp(G)$ & $8$ & $4$ & $8$ & \!\!$4$ or $8$\!\! & $8$ & $8$ & $8$ & $16$ & $8$ & $8$ & $4$\\\hline
$B_0(G)$ & \multicolumn{10}{c|}{$\bZ/2\bZ$} & $(\bZ/2\bZ)^{\oplus 2}$\!\\\hline 
\# of $G$'s & 48 & 55 & 18 & 6 & 26 & 20 & 10 & 9 & 2 & 2 & 34\\\hline
\end{tabular}
\end{table}
%


It is natural to ask the (stably) birational classification 
of $\bC(G)$ for groups $G$ of order $128$. 
In particular, 
what happens to $\bC(G)$ with $B_0(G)\neq 0$? 
The following theorem (Theorem \ref{thm}) 
gives a partial answer to this question. 

\begin{definition}\label{defL23}
Let $k$ be a field with {\rm char} $k\neq 2$ and 
$k(X_1,X_2,X_3,X_4,X_5,X_6,$ $X_7)$ be the rational function field 
over $k$ with variables $X_1,X_2,X_3,X_4,X_5,X_6,X_7$.\\
{\rm (1)} {\it The field $L_k^{(2)}$} is defined to be 
$k(X_1,X_2,X_3,X_4,X_5,X_6)^{\langle\rho\rangle}$ 
where $\langle\rho\rangle\simeq C_4$ acts on\\ 
$k(X_1,X_2,X_3,X_4,X_5,X_6)$ by $k$-automorphisms 
\begin{align*}
\rho : &\ X_1 \mapsto X_2, X_2 \mapsto -X_1, X_3 \mapsto X_4, 
X_4 \mapsto X_3,\\
&\ X_5 \mapsto X_6, 
 X_6 \mapsto \frac{(X_1^2X_2^2-1)(X_1^2X_3^2+X_2^2-X_3^2-1)}{X_5}.
\end{align*}
{\rm (2)} {\it The field $L_k^{(3)}$} 
is defined to be $k(X_1,X_2,X_3,X_4,X_5,X_6,X_7)^{\langle\lambda_1,\lambda_2\rangle}$ 
where $\langle\lambda_1,\lambda_2\rangle \simeq C_2\times C_2$ acts on 
$k(X_1,X_2,X_3,X_4,X_5,X_6,X_7)$ by $k$-automorphisms 
\begin{align*}
\lambda_1 : &\ X_1 \mapsto X_1, X_2 \mapsto \frac{X_1}{X_2}, 
X_3 \mapsto \frac{1}{X_1 X_3}, X_4 \mapsto \frac{X_2 X_4}{X_1 X_3},\\ 
&\ X_5 \mapsto -\frac{X_1 X_6^2-1}{X_5}, X_6 \mapsto -X_6, X_7 \mapsto X_7,\\
\lambda_2 : &\ X_1 \mapsto \frac{1}{X_1}, X_2 \mapsto X_3, X_3 \mapsto X_2, 
 X_4 \mapsto \frac{(X_1 X_6^2-1) (X_1 X_7^2-1)}{X_4},\\
&\ X_5 \mapsto -X_5, X_6 \mapsto -X_1 X_6, X_7 \mapsto -X_1 X_7.
\end{align*}
\end{definition}

\begin{theorem}[Hoshi {\cite[Theorem 1.31]{Hos16}}]\label{thm}
Let $G$ be a group of order $128$. 
Assume that $B_0(G)\neq 0$. 
Then $\bC(G)$ and $L_\bC^{(m)}$ are stably $\bC$-isomorphic where 
\begin{align*}
m=
\begin{cases}
1&{\rm if}\ G\ {\rm belongs\ to\ } \Phi_{16}, \Phi_{31}, \Phi_{37}, 
\Phi_{39}, \Phi_{43}, \Phi_{58}, \Phi_{60}\ {\rm or}\ \Phi_{80},\\
2&{\rm if}\ G\ {\rm belongs\ to\ } \Phi_{106}\ {\rm or}\ \Phi_{114},\\
3&{\rm if}\ G\ {\rm belongs\ to\ } \Phi_{30}.
\end{cases}
\end{align*}
In particular, ${\rm Br}_{\rm nr}(L_\bC^{(1)})\simeq {\rm Br}_{\rm nr}(L_\bC^{(2)})\simeq\bZ/2\bZ$ 
and ${\rm Br}_{\rm nr}(L_\bC^{(3)})\simeq (\bZ/2\bZ)^{\oplus 2}$ 
and hence the fields 
$L_\bC^{(1)}$, $L_\bC^{(2)}$ and $L_\bC^{(3)}$ are not $($retract, stably$)$ $\bC$-rational.
\end{theorem}

For $m=1,2$, 
the fields $L_\bC^{(m)}$ and $L_\bC^{(3)}$ are not stably $\bC$-isomorphic 
because their unramified Brauer groups are not isomorphic. 
However, we do not know whether the fields $L_\bC^{(1)}$ and $L_\bC^{(2)}$ are stably 
$\bC$-isomorphic. 
If not, it is interesting to evaluate the higher unramified cohomologies. 

\section{Rationality problem for multiplicative invariant fields}

Let $k$ be a field, $G$ be a finite group and
$\rho : G\rightarrow GL(V)$ be a faithful representation of $G$ 
where $V$ is a finite-dimensional vector space over $k$.
Then $G$ acts on the rational function field $k(V)$.

We consider the rationality problem for $k(V)^G$. 
By No-name Lemma (cf. Miyata \cite[Remark 3]{Miy71}), 
it is known that $k(G)$ is stably $k$-rational if and only if so 
is $k(V)^G$ where $\rho : G\rightarrow GL(V)$ 
is any faithful representation of $G$ over $k$. 
Thus the rationality problem of $k(V)^G$ over $k$ is 
also called Noether's problem. 

In order to solve the rationality problem of $k(V)^G$, 
it is natural and almost inevitable that we reduce the problem 
to that of the multiplicative invariant field $k(M)^G$ defined 
in Definition \ref{d1.2}; 
an illustration of reducing Noether's problem to the multiplicative 
invariant field can be found in, 
e.g. \cite{CHKK10}, \cite[Example 13.7]{HKY11}. 

When $M$ is a $G$-lattice with ${\rm rank}_\bm{Z} M=n$, 
the multiplicative invariant field $k(M)^G$ is nothing but 
$k(x_1,\ldots, x_n)^G$, the fixed field of the rational function 
field $k(x_1,\ldots, x_n)$ on which $G$ acts by multiplicative actions. 
\begin{definition}\label{d1.1}
Let $G$ be a finite group and $\bm{Z}[G]$ be the group ring. 
A finitely generated $\bm{Z}[G]$-module $M$ is called a 
{\it $G$-lattice} if, as an abelian group, 
$M$ is a free abelian group of finite rank. 
We will write ${\rm rank}_\bm{Z} M$ for the rank of $M$ 
as a free abelian group.
A $G$-lattice $M$ is called {\it faithful} if, for any 
$\sigma\in G\setminus\{1\}$, $\sigma\cdot x\neq x$ for some $x\in M$.

Suppose that $G$ is any finite group and $\Phi : G\rightarrow GL_n(\bm{Z})$ is a group homomorphism, i.e. an integral representation of $G$. Then the group $\Phi(G)$ acts naturally on the free abelian group $M:=\bm{Z}^{\oplus n}$; thus $M$ becomes a $\bm{Z}[G]$-module. We call $M$ the 
{\it $G$-lattice associated to $\Phi$ $($or $\Phi(G)$$)$}. 
Conversely, if $M$ is a $G$-lattice with ${\rm rank}_\bm{Z} M=n$, write $M=\oplus_{1\leq i\leq n}\bm{Z}\cdot x_i$.
Then there is a group homomorphism
$\Phi : G\rightarrow GL_n(\bm{Z})$ defined as follows: If
$\sigma\cdot x_i=\sum_{1\leq j\leq n}a_{ij}\,x_j$
where $\sigma\in G$ and $a_{ij}\in\bm{Z}$,
define $\Phi(\sigma)=(a_{ij})_{1\leq i,j\leq n}\in GL_n(\bm{Z})$.

When the group homomorphism $\Phi : G\rightarrow GL_n(\bm{Z})$ is injective, the corresponding $G$-lattice is a faithful $G$-lattice. For examples, any finite subgroup $G$ of $GL_n(\bm{Z})$ gives rise to
a faithful $G$-lattice of rank $n$.

The list of all the finite subgroups of $GL_n(\bm{Z})$ (with $n \le 4$), up to conjugation, can be found in the book \cite{BBNWZ78} and in GAP \cite{GAP}. 
As to the situations of $GL_n(\bm{Z})$ (with $n \ge 5$), Plesken etc. found the lists of all the finite subgroups of $GL_n(\bm{Z})$ (with $n =5$ and $6$); see \cite{PS00} and the references therein. These lists may be found in the GAP 
package CARAT \cite{CARAT} and also in \cite[Chapter 3]{HY17}.

Here is a list of the total number of lattices,
up to isomorphism, of a given rank:\vspace*{2mm} 
\begin{center}
\begin{tabular}{|c|cccccc|}\hline
rank & $1$ & $2$ & $3$ & $4$ & $5$ & $6$\\\hline
\# of $G$-lattices & $2$ & $13$ & $73$ & $710$ & $6079$ & $85308$\\\hline
\end{tabular}\vspace*{2mm}
\end{center}
\end{definition}

\begin{definition}\label{d1.2}
Let $M$ be a $G$-lattice of rank $n$ and write
$M=\oplus_{1\leq i\leq n}\bm{Z}\cdot x_i$.
For any field $k$, define $k(M)=k(x_1,\ldots,x_n)$
the rational function field over $k$ with $n$ variables $x_1,\ldots,x_n$.
Define a {\it multiplicative action} of $G$ on $k(M)$: 
For any $\sigma\in G$,
if $\sigma\cdot x_i=\sum_{1\leq j\leq n}a_{ij}\,x_j$ in the $G$-lattice $M$,
then we define $\sigma\cdot x_i=\prod_{1\leq j\leq n}x_j^{a_{ij}}$
in the field $k(M)$.
Note that $G$ acts trivially on $k$.
The above multiplicative action is called a {\it purely monomial action} of $G$
on $k(M)$ in \cite{HK92} 
and $k(M)^G$ is called a {\it multiplicative invariant field} 
in \cite{Sal87}.


When $M$ is the $G$-lattice $\bm{Z}[G]$ 
where $M=\oplus_{g \in G} \bm{Z} \cdot x_g$ and 
$h \cdot x_g=x_{hg}$ for $h, g \in G$, we have 
$k(M)=k(x_g \mid g \in G)$ and 
$k(M)^G=k(G)$ (see Section 1). 
Note that $k(G)= k(V_{\rm reg})^G$ where $G \rightarrow GL(V_{\rm reg})$ is the regular representation of $G$ over $k$.
\end{definition}
%
\begin{theorem}[{Hajja \cite{Haj87}}]\label{thHaj87}
Let $k$ be a field and $G$ be a finite group acting on $k(x_1,x_2)$
by monomial $k$-automorphisms.
Then $k(x_1,x_2)^{G}$ is $k$-rational.
\end{theorem}
\begin{theorem}[{Hajja and Kang \cite{HK92, HK94}, Hoshi and Rikuna \cite{HR08}}]
\label{thHKHR}
Let $k$ be a field and $G$ be a finite group acting on
$k(x_1,x_2,x_3)$ by purely monomial $k$-automorphisms.
Then $k(x_1,x_2,x_3)^G$ is $k$-rational.
\end{theorem}
\begin{theorem}[Hoshi, Kang and Kitayama {\cite[Theorem 1.16]{HKK14}}]\label{th116}
Let $k$ be a field, $G$ be a finite group and $M$ be a $G$-lattice
with ${\rm rank}_{\bZ} M=4$ such that $G$ acts on $k(M)$ by
purely monomial $k$-automorphisms.
If $M$ is decomposable, i.e. $M=M_1\oplus M_2$ as $\bZ[G]$-modules where
$1\le {\rm rank}_{\bZ} M_1 \le 3$, then $k(M)^G$ is $k$-rational.
\end{theorem}
\begin{theorem}[Hoshi, Kang and Kitayama {\cite[Theorem 6.2]{HKK14}}]\label{th62}
Let $k$ be a field, $G$ be a finite group and $M$ be a $G$-lattice
such that $G$ acts on $k(M)$ by purely monomial $k$-automorphisms.
Assume that
{\rm (i)} $M=M_1\oplus M_2$ as $\bZ[G]$-modules where
${\rm rank}_{\bZ}M_1=3$ and ${\rm rank}_{\bZ} M_2=2$,
{\rm (ii)} either $M_1$ or $M_2$ is a faithful $G$-lattice.
Then $k(M)^G$ is $k$-rational except the following situation:
{\rm char} $k\ne 2$, $G=\langle\sigma,\tau\rangle \simeq D_4$
and $M_1=\bigoplus_{1\le i\le 3}
\bZ x_i$, $M_2=\bigoplus_{1\le j\le 2} \bZ y_j$ such that
$\sigma:x_1\leftrightarrow x_2$, $x_3\mapsto -x_1-x_2-x_3$,
$y_1\mapsto y_2\mapsto -y_1$, $\tau: x_1\leftrightarrow x_3$,
$x_2\mapsto -x_1-x_2-x_3$, $y_1\leftrightarrow y_2$ where the
$\bZ[G]$-module structure of $M$ is written additively.
For the exceptional case, $k(M)^G$ is not retract $k$-rational.
\end{theorem}
%
\begin{definition}\label{d1.3}
Let $k$ be a field and $\mu$ be a multiplicative subgroup of $k \setminus \{0 \}$ containing all the
roots of unity in $k$.
If $M$ is a $G$-lattice, a {\it $\mu$-extension} is an exact sequence of
$\bm{Z}[G]$-modules given by
$(\alpha) : 1\rightarrow\mu\rightarrow M_\alpha\rightarrow M\rightarrow 0$
where $G$ acts trivially on $\mu$. Be aware that $M_{\alpha}=\mu \oplus M$ as abelian groups, but not as $\bm{Z}[G]$-modules except when the extension $(\alpha)$ splits.

As in Definition \ref{d1.2}, if $M=\oplus_{1\leq i\leq n}\bm{Z}\cdot x_i$ and $M_{\alpha}$ is a $\mu$-extension,
we define the field $k_\alpha(M)=k(x_1,\ldots,x_n)$ the rational function
field over $k$ with $n$ variables $x_1,\ldots,x_n$; 
the action of $G$ on $k_\alpha(M)$ will be described in the next paragraph. Note that $M_{\alpha}$ is embedded into the
multiplicative group $k_\alpha(M) \setminus \{ 0 \}$ by sending $(\epsilon, \sum_{1 \le i \le n}b_ix_i)\in \mu \oplus M$ to the element $\epsilon \prod_{1 \le i \le n}x_i^{b_i}$ in the field $k_\alpha(M) =k(x_1,\ldots,x_n)$.

The group $G$ acts on $k_\alpha(M)$ by a twisted multiplicative action:
Suppose that, in $M$ we have $\sigma\cdot x_i=\sum_{1\leq j\leq n}a_{ij}\,x_j$, and
in $M_\alpha$ we have $\sigma\cdot x_i=\varepsilon_i(\sigma)+\sum_{1\leq j\leq n}a_{ij}\,x_j$ where $\varepsilon_i(\sigma)\in\mu$.
Then we define
$\sigma\cdot x_i=\varepsilon_i(\sigma)\prod_{1\leq j\leq n}x_j^{a_{ij}}$
in $k_\alpha(M)$. Again $G$ acts trivially on the coefficient field $k$.
The above group action is called {\it monomial group action} in \cite{HK92} 
and $k_\alpha(M)^G$ is called {\it twisted multiplicative 
invariant field} in \cite{Sal90}.

Note that, if the extension $(\alpha) : 1\rightarrow \mu\rightarrow M_\alpha
\rightarrow M\rightarrow 0$ is a split extension, then $k_\alpha(M)=k(M)$
and the twisted multiplicative action is reduced to the multiplicative
action in Definition \ref{d1.2}.
\end{definition}


For any faithful linear representation $G \to GL(V)$ 
of $G$, 
we have 
${\rm Br}_{\rm nr}(\bC(V)^G) \simeq B_0(G)$ 
by No-name Lemma (see \cite{Sal90}). 

The formula in \cite[Theorem 12]{Sal90} 
(Theorem \ref{thBog})
can be used to compute not only ${\rm Br}_{\rm nr}(\bm{C}(V)^G)$, but also ${\rm Br}_{\rm nr}(\bm{C}_\alpha(M)^G)$ where $\bm{C}_\alpha(M)$ is the rational function field associated to the $\mu$-extension $M_{\alpha}$: 
\begin{theorem}[Saltman {\cite[Theorem 12]{Sal90}}]\label{thSa4}
Let $k$ be an algebraically closed field with char $k =0$, and $G$ be a finite group. If $M$ is a $G$-lattice and $(\alpha): 1 \to \mu \to M_{\alpha} \to M \to 0$ is a $\mu$-extension such that {\rm (i)} $M$ is a faithful $G$-lattice, and {\rm (ii)} $H^2(G,\mu) \to H^2(G,M_\alpha)$ is injective, then 
\[
{\rm Br}_{\rm nr}(k_\alpha(M)^G)=\bigcap_{A} {\rm Ker}\{
{\rm res} : H^2(G,M_\alpha)\rightarrow H^2(A,M_\alpha)\}
\]
where $A$ runs over all the bicyclic subgroups of $G$.

In particular, if the $\mu$-extension $(\alpha) : 1\rightarrow\mu\rightarrow M_\alpha\rightarrow M\rightarrow 0$ splits, then ${\rm Br}_{\rm nr}(k(M)^G)\simeq B_0(G) \oplus \bigcap_{A} {\rm Ker}\{
{\rm res} : H^2(G,M)\rightarrow H^2(A,M)\}$ where $A$ runs over bicyclic subgroups of $G$.
\end{theorem}
\begin{definition}\label{d1.4}
By Definition \ref{defSal}, 
${\rm Br}_{\rm nr}(K)$ is a subgroup of the Brauer group ${\rm Br}(K)$. On the other hand, 
the map of the Brauer groups ${\rm Br}(k_\alpha(M)^G) \to {\rm Br}(k_\alpha(M))$ sends ${\rm Br}_{\rm nr}(k_\alpha(M)^G)$ to ${\rm Br}_{\rm nr}(k_\alpha(M))$ \cite[Theorem 2.1]{Sal87}. Since ${\rm Br}_{\rm nr}(k_\alpha(M))=0$ by \cite[Proposition, 2.2]{Sal87}, it follows that the unramified Brauer group ${\rm Br}_{\rm nr}(k_\alpha(M)^G)$ is a subgroup of the relative Brauer group ${\rm Br}(k_\alpha(M)/k_\alpha(M)^G)$. 
As ${\rm Br}(k_\alpha(M)/k_\alpha(M)^G)$ is isomorphic to the cohomology group $H^2(G,k_\alpha(M)^{\times})$, we may regard ${\rm Br}_{\rm nr}(k_\alpha(M)^G)$ as a subgroup of $H^2(G,k_\alpha(M)^{\times})$.

Through the embedding $M_\alpha\hookrightarrow k_\alpha(M)^\times$,
there is a canonical injection $H^2(G,M_\alpha)$ $\hookrightarrow$ ${\rm Br}(k_\alpha(M)^G)$ \cite[page 536]{Sal90}. 
Identifying 
${\rm Br}_{\rm nr}(k_\alpha(M)^G)$ and $H^2(G,M_\alpha)$ as subgroups of $H^2(G,k_\alpha(M)^{\times})$, 
we see that 
${\rm Br}_{\rm nr}(k_\alpha(M)^G)$ is a subgroup of $H^2(G,M_\alpha)$ \cite[page 536]{Sal90}. Thus we write $H^2_{\rm nr}(G,M_\alpha)$ for ${\rm Br}_{\rm nr}(k_\alpha(M)^G)$ (see \cite{Sal90}).

Note that there is a natural map $H^2(G,\bm{Q}/\bm{Z}) \to H^2(G,M_\alpha)$.
Clearly this map is injective if the $\mu$-extension $(\alpha) : 1\rightarrow\mu\rightarrow M_\alpha\rightarrow M\rightarrow 0$ splits. In this case, regarding $H^2(G,\bm{Q}/\bm{Z})$ and $H^2(G,M)$ as subgroups of $H^2(G,M_\alpha)$, we define $H^2_{\rm nr}(G,\bm{Q}/\bm{Z})=H^2(G,\bm{Q}/\bm{Z}) \cap
{\rm Br}_{\rm nr}(k_\alpha(M)^G)$ and $H^2_{\rm nr}(G,M)=H^2(G,M) \cap
{\rm Br}_{\rm nr}(k_\alpha(M)^G)$.
It follows that
${\rm Br}_{\rm nr}(k_\alpha(M)^G)=H_{\rm nr}^2(G,\bQ/\bZ)\oplus H_{\rm nr}^2(G,M)$.
By Theorems \ref{thBog} and \ref{thSa4},
we have $H^2_{\rm nr}(G,\bQ/\bZ)\simeq B_0(G)$ and
$H^2_{\rm nr}(G,M)\simeq \bigcap_{A} {\rm Ker}\{
{\rm res} : H^2(G,M)\rightarrow H^2(A,M)\}$ where $A$ runs over
bicyclic subgroups of $G$.
\end{definition}
\begin{theorem}[{Barge \cite[Theorem II.7]{Bar89}}]\label{thB1}
Let $G$ be a finite group.
Then the following conditions are equivalent:\\
{\rm (1)} All the Sylow subgroups of $G$ are bicyclic;\\
{\rm (2)} ${\rm Br}_{\rm nr}(\bm{C}(M)^G)=0$ for all $G$-lattices $M$.
\end{theorem}
\begin{theorem}[Barge {\cite[Theorem IV-1]{Bar97}}] \label{thB2}
Let $G$ be a finite group.
Then the following conditions are equivalent:\\
{\rm (1)} All the Sylow subgroups of $G$ are cyclic;\\
{\rm (2)} ${\rm Br}_{\rm nr}(\bm{C}_{\alpha}(M)^G)=0$ for all $G$-lattices $M$,
for all short exact sequences of
$\bm{Z}[G]$-modules $\alpha : 0 \rightarrow \bm{C}^{\times}
\rightarrow M_{\alpha} \rightarrow M \rightarrow 0$.
\end{theorem}
As in Definition \ref{d1.4}, we have 
${\rm Br}_{\rm nr}(\bm{C}(M)^G) \simeq B_0(G) \oplus H^2_{\rm nr}(G,M)$ where $B_0(G)$ is the Bogomolov multiplier 
and $H^2_{\rm nr}(G,M)\leq H^2(G,M)$. We remark that $B_0(G)$ is related to the rationality of $\bC(V)^G$ where $G \to GL(V)$ is any faithful linear representation of $G$ over $\bC$; on the other hand, $H^2_{\rm nr}(G,M)$ arises from the multiplicative nature of the field $\bm{C}(M)^G$. 

In case ${\rm rank}_\bm{Z} M \le 3$, ${\rm Br}_{\rm nr}(\bm{C}(M)^G)=0$ 
for all $G$-lattices $M$ because $\bm{C}(M)^G$ are always $\bC$-rational 
(see Theorem \ref{thHaj87} and Theorem \ref{thHKHR}). 
The following theorem \cite[Theorem 1.10]{HKY} gives the classification of 
all the lattices $M$ with ${\rm Br}_{\rm nr}(\bm{C}(M)^G)\neq 0$ 
when ${\rm rank}_\bm{Z} M\leq 6$. 
Thus $\bm{C}(M)^G$ are not retract $\bC$-rational 
for these lattices (and thus are not $\bC$-rational). 

Let $C_n$ (resp. $D_n$, $QD_{8n}$, $Q_{8n}$) 
be the cyclic group of order $n$ 
(resp. the dihedral group of order $2n$, 
the quasi-dihedral group of order $16n$, 
the generalized quaternion group of order $8n$). 
\begin{theorem}[{Hoshi, Kang and Yamasaki \cite[Theorem 1.10]{HKY}}]\label{t1.5}
Let $G$ be a finite group and $M$ be a faithful $G$-lattice.
\\
{\rm (1)} If ${\rm rank}_\bm{Z} M \le 3$, then ${\rm Br}_{\rm nr}(\bm{C}(M)^G)=0$.\\
{\rm (2)} If ${\rm rank}_\bm{Z} M =4$, then ${\rm Br}_{\rm nr}(\bm{C}(M)^G)\neq 0$
if and only if $M$ is one of the $5$ cases in {\rm Table} $1$. 
Moreover, if $M$ is one of the $5$ $G$-lattices with 
${\rm Br}_{\rm nr}(\bm{C}(M)^G)\neq 0$, then $B_0(G)=0$ and ${\rm Br}_{\rm nr}(\bm{C}(M)^G)=H_{\rm nr}^2(G,M)$.\\
{\rm (3)} If ${\rm rank}_\bm{Z} M =5$, then ${\rm Br}_{\rm nr}(\bm{C}(M)^G)\neq 0$ if and only if
$M$ is one of the $46$ cases in {\rm \cite[Table $2$]{HKY}}.
Moreover, if $M$ is one of the $46$ $G$-lattices with 
${\rm Br}_{\rm nr}(\bm{C}(M)^G)\neq 0$, then $B_0(G)=0$ and ${\rm Br}_{\rm nr}(\bm{C}(M)^G)=H_{\rm nr}^2(G,M)$. \\
{\rm (4)} If ${\rm rank}_\bm{Z} M =6$, then
${\rm Br}_{\rm nr}(\bm{C}(M)^G)\neq 0$ if and only if
$M$ is one of the $1073$ cases as in {\rm \cite[Table $3$]{HKY}}. 
Moreover, if $M$ is one of the 
$1073$ $G$-lattices with 
${\rm Br}_{\rm nr}(\bm{C}(M)^G)\neq 0$, then $B_0(G)=0$ and ${\rm Br}_{\rm nr}(\bm{C}(M)^G)=H_{\rm nr}^2(G,M)$,
except for $24$ cases
with $B_0(G)=\bm{Z}/2\bm{Z}$ where the {\rm CARAT ID} of $G$ are
$(6,6458,i)$, $(6,6459,i)$, $(6,6464,i)$ $(1\leq i\leq 8)$.
Note that $22$ cases out of the exceptional $24$ cases satisfy 
$H_{\rm nr}^2(G,M)=0$.
\end{theorem}

\begin{center}
Table $1$: $5$ $G$-lattices $M$ of rank $4$ with 
 ${\rm Br}_{\rm nr}(\bm{C}(M)^G)\neq 0$\vspace*{2mm}\\
\begin{tabular}{|lllcl|}
\hline
$G(n,i)$ & $G$ & GAP ID & $B_0(G)$ & $H_{\rm nr}^2(G,M)$\\\hline
$(8,3)$ & $D_4$ & $(4,12,4,12)$ & $0$ & $\bZ/2\bZ$\\
$(8,4)$ & $Q_8$ & $(4,32,1,2)$ & $0$ & $(\bZ/2\bZ)^{\oplus 2}$\\
$(16,8)$ & $QD_8$ & $(4,32,3,2)$ & $0$ & $\bZ/2\bZ$\\
$(24,3)$ & $SL_2(\bm{F}_3)$ & $(4,33,3,1)$ & $0$ & $(\bZ/2\bZ)^{\oplus 2}$\\
$(48,29)$ & $GL_2(\bm{F}_3)$ & $(4,33,6,1)$ & $0$ & $\bZ/2\bZ$\\\hline
\end{tabular}\vspace*{2mm}
\end{center}

\begin{remark}\label{remm}
(1) The above theorem remains valid if we replace the coefficient field $\bm{C}$ by any algebraically closed field $k$ with char $k=0$.\\
(2) If $M$ is of rank $\leq 6$ and ${\rm Br}_{\rm nr}(\bm{C}(M^G))\neq 0$,
then $G$ is solvable and non-abelian, and
${\rm Br}_{\rm nr}(\bm{C}(M)^G)\simeq \bZ/2\bZ$, $\bZ/3\bZ$ or
$\bZ/2\bZ\oplus\bZ/2\bZ$. 
The case where 
${\rm Br}_{\rm nr}(\bm{C}(M)^G)\simeq\bZ/3\bZ$ occurs only for
$4$ groups $G$ of order $27$, $27$, $54$, $54$
with the {\rm CARAT ID}
$(6,2865,1)$, $(6,2865,3)$, $(6,2899,3)$, $(6,2899,5)$
which are isomorphic to $C_9\rtimes C_3$, $C_9\rtimes C_3$,
$(C_9\rtimes C_3)\rtimes C_2$, $(C_9\rtimes C_3)\rtimes C_2$ respectively.
For {\rm CARAT ID}, see Hoshi and Yamasaki \cite[Chapter 3]{HY17}.
\\
(3)
The group $G$ $(\simeq D_4)$ which appears as the exceptional case in Theorem \ref{th62} (i.e. \cite[Theorem 6.2]{HKK14}) satisfies the property that ${\rm Br}_{\rm nr}(\bm{C}(M)^G)= H^2_{\rm nr} (G,M)\neq 0$ where $M$ is the associated lattice. It follows that $\bC(M)^G$ is not retract rational.

In Theorem \ref{th62}, note that both $\bC(M_1)^G$ and $\bC(M_2)^G$ are rational by Theorem \ref{thHKHR} and Theorem \ref{thHaj87}. Thus
${\rm Br}_{\rm nr}(\bm{C}(M_2)^G)=0$ and 
$H^2 _{\rm nr} (G,M_2)=0$. But $M_1$ is not a faithful $G$-lattice and we cannot apply Theorem \ref{thSa4} to $\bC(M_1)^G$. Hence it is possible that $H^2 _{\rm nr} (G,M_1)$ is non-trivial. Because $H^2 _{\rm nr} (G,M) \simeq H^2 _{\rm nr} (G,M_1)\oplus H^2 _{\rm nr} (G,M_2)$, this allows for the possibility that $H^2 _{\rm nr} (G,M)$ is non-trivial. 
Indeed, 
it can be shown that $H^2 _{\rm nr} (G,M_1) \simeq \bZ/2 \bZ$ and therefore ${\rm Br}_{\rm nr}(\bm{C}(M)^G)= H^2 _{\rm nr} (G,M_1)\simeq \bZ/2 \bZ$.\\ 
{\rm (4)} Here is a summary of Theorem \ref{t1.5}:
\begin{center}
\begin{tabular}{|l|rrrrrr|}\hline
${\rm rank}_\bm{Z} M$ & $1$ & $2$ & $3$ & $4$ & $5$ & $6$\\\hline
\# of $G$-lattices $M$ & $2$ & $13$ & $73$ & $710$ & $6079$ & $85308$\\\hline
\# of $G$-lattices $M$ with ${\rm Br}_{\rm nr}(\bm{C}(M)^G)\neq 0$
& $0$ & $0$ & $0$ & $5$ & $46$ & $1073$\\\hline
\end{tabular}\vspace*{2mm}
\end{center}
\end{remark}
\begin{theorem}[{Hoshi, Kang and Yamasaki \cite[Theorem 4.4]{HKY}}]\label{equiv} 
The following fields $K$ are stably equivalent each other:\\
{\rm (1)} $\bm{C}(G)$ where $G$ is a group of order $64$ which belongs to the $16$th isoclinism class $\Phi_{16}$ $($see the $9$ groups defined as in 
{\rm Theorem \ref{thCHKK10} (1)}$)$;\\
{\rm (2)} $\bm{C}(x_1,x_2,x_3,x_4)^{D_4}$
where $D_4=\langle\sigma,\tau\rangle$ acts on
$\bm{C}(x_1,x_2,x_3,x_4)$ by
\begin{align*}
&\sigma: x_1\mapsto x_2x_3, x_2\mapsto x_1x_3, x_3\mapsto x_4, x_4\mapsto \tfrac{1}{x_3},\\
&\tau: x_1\mapsto \tfrac{1}{x_2}, x_2\mapsto\tfrac{1}{x_1},
x_3\mapsto\tfrac{1}{x_4}, x_4\mapsto\tfrac{1}{x_3}
\end{align*}
$($see {\rm Theorem \ref{t1.5} (2) and Table $1$}$)$;\\
{\rm (3)} $\bm{C}(y_1,y_2,y_3,y_4,y_5)^{D_4}$ where $D_4=\langle\sigma,\tau\rangle$ acts on $\bm{C}(y_1,y_2,y_3,y_4,y_5)$ by
\begin{align*}
&\sigma: y_1\mapsto y_2, y_2\mapsto y_1, y_3\mapsto \tfrac{1}{y_1y_2y_3},
y_4\mapsto y_5, y_5\mapsto \tfrac{1}{y_4},\\
&\tau: y_1\mapsto y_3, y_2\mapsto \tfrac{1}{y_1y_2y_3}, y_3\mapsto y_1,
y_4\mapsto y_5, y_5\mapsto y_4
\end{align*}
$($see {\rm Theorem \ref{th62}}$)$;\\
{\rm (4)} $\bm{C}(z_1,z_2,z_3,z_4)^{C_2\times C_2}$ where
$C_2\times C_2=\langle\sigma,\tau\rangle$ acts on
$\bm{C}(z_1,z_2,z_3,z_4)$ by
\begin{align*}
&\sigma: z_1\mapsto z_2, z_2\mapsto z_1, z_3\mapsto \tfrac{1}{z_1z_2z_3},
z_4\mapsto \tfrac{-1}{z_4},\\
&\tau: z_1\mapsto z_3, z_2\mapsto \tfrac{1}{z_1z_2z_3}, z_3\mapsto z_1,
z_4\mapsto -z_4
\end{align*}
$($see {\rm \cite[Proof of Theorem 6.4]{HKK14}}$)$;\\
{\rm (5)} $\bm{C}(w_1,w_2,w_3,w_4)^{C_2}$ where
$C_2=\langle\sigma\rangle$ acts on
$\bm{C}(w_1,w_2,w_3,w_4)$ by
\begin{align*}
&\sigma: w_1\mapsto -w_1, w_2\mapsto \tfrac{w_4}{w_2},
w_3\mapsto \tfrac{(w_4-1)(w_4-w_1^2)}{w_3}, w_4\mapsto w_4
\end{align*}
$($see {\rm \cite[Theorem 6.3]{HKK14}}$)$. 

In particular, the unramified cohomology groups
$H_{\rm nr}^i(K,\bm{Q}/\bm{Z})$ of the fields $K$ in {\rm (1)}--{\rm (5)} 
coincide and ${\rm Br}_{\rm nr}(K)\simeq \bm{Z}/2\bm{Z}$.
\end{theorem}

As in Remark \ref{remm} (2), 
all the $G$-lattices $M$ with ${\rm rank}_\bZ M\leq 6$ and $H_{\rm nr}^2(G,M)\neq 0$ in Theorem \ref{t1.5} 
satisfy the condition that $G$ is non-abelian and solvable. 
Examples of $G$-lattices $M$ with $H_{\rm nr}^2(G,M)\neq 0$ 
where $G$ is abelian 
(resp. non-solvable; in fact, simple) 
are given in \cite{HKY} as follows: 

\begin{theorem}[{Hoshi, Kang and Yamasaki \cite[Theorem 6.1]{HKY}}]\label{t5.4}
Let $G$ be an elementary abelian group of order $2^n$ in $GL_7(\bm{Z})$ and
$M$ be the associated $G$-lattice of rank $7$.
Then 
${\rm Br}_{\rm nr}(\bm{C}(M)^G)\neq 0$ if and only if $G$ is isomorphic up to conjugation to one of the nine groups $G_1,\ldots,G_9\leq GL_7(\bm{Z})$ 
as in {\rm \cite[Theorem 6.1]{HKY}} 
where each of $G_i$ is isomorphic to $(C_2)^3$ as an abstract group. 
Moreover, ${\rm Br}_{\rm nr}(\bm{C}(M)^{G_i})=H_{\rm nr}^2(G_i,M)\simeq \bZ/2\bZ$
$($resp. $\bZ/2\bZ\oplus\bZ/2\bZ$$)$ for $1\leq i\leq 8$ $($resp. $i=9$$)$.
\end{theorem}

\begin{theorem}[{Hoshi, Kang and Yamasaki \cite[Theorem 6.2]{HKY}}]\label{t5.5}
Embed $A_6$\\
into $S_{10}$ through the isomorphism $A_6\simeq PSL_2(\bF_9)$, which acts on the projective line $\bm{P}^1 _{\bF_9}$ via fractional linear transformations. Thus we may regard $A_6$ as a transitive subgroup of $S_{10}$.  Let $N=\oplus_{1\leq i\leq 10}\bZ\cdot x_i$ be the  $S_{10}$-lattice defined by $\sigma \cdot x_i = x_{\sigma(i)}$ for any $\sigma \in S_{10}$; it becomes an $A_6$-lattice by restricting the action of $S_{10}$ to $A_6$. Define $M=N/(\bZ\cdot \sum_{i=1}^{10}x_i)$ with ${\rm rank}_\bZ M=9$.
There exist exactly six $A_6$-lattices $M=M_1$, $M_2,\ldots,M_6$
which are $\bQ$-conjugate but not $\bZ$-conjugate to each other; in fact, all these $M_i$ form a single $\bQ$-class, but this $\bQ$-class consists of six $\bZ$-classes.
Then we have
\begin{align*}
H_{\rm nr}^2(A_6,M_1)\simeq H_{\rm nr}^2(A_6,M_3)\simeq\bZ/2\bZ,\quad
H_{\rm nr}^2(A_6,M_i)=0\ {\rm for}\ i=2,4,5,6.
\end{align*}

In particular, $\bC(M_1)^{A_6}$ and $\bC(M_3)^{A_6}$
are not retract $\bC$-rational.
Furthermore, 
the lattices $M_1$ and $M_3$ may be distinguished by the 
Tate cohomology groups: 
\begin{align*}
&H^1(A_6,M_1)=0,& &\hspace*{-1cm}\widehat{H}^{-1}(A_6,M_1)=\bZ/10\bZ,\\ 
&H^1(A_6,M_3)=\bZ/5\bZ,& &\hspace*{-1cm}\widehat{H}^{-1}(A_6,M_3)=\bZ/2\bZ.
\end{align*}
\end{theorem}
Motivated by the $G$-lattices in Theorem \ref{t1.5} (2) (see Table $1$), 
the following $G$-lattices $M$ of rank $2n+2$, $4n$ and $p(p-1)$ 
($n$ is any positive integer and $p$ is any odd prime number) 
with ${\rm Br}_{\rm nr}(\bC(M)^G)\neq 0$ were constructed in \cite{HKY}: 
%
\begin{theorem}[{Hoshi, Kang and Yamasaki \cite[Theorem 7.2]{HKY}}]\label{thmD4n}
Let $G$ $=$ $\langle \sigma$,\\ 
$\tau \mid \sigma^{4n}=\tau^2=1,
\tau^{-1}\sigma\tau=\sigma^{-1}\rangle\simeq D_{4n}$,
the dihedral group of order $8n$ where $n$ is any positive integer.
Let $M$ be the $G$-lattice of rank $2n+2$ defined in 
{\rm \cite[Definition 7.1]{HKY}}.
Then $H_{\rm nr}^2(G,M)\simeq \bZ/2\bZ$.
Consequently, $\bC(M)^G$ is not retract $\bC$-rational.
\end{theorem}
\begin{theorem}[{Hoshi, Kang and Yamasaki \cite[Theorem 7.5]{HKY}}]\label{thmQDQ}\\
{\rm (1)} Let $n$ be any positive integer and 
$G=\langle \sigma,\tau \mid \sigma^{8n}=\tau^2=1,
\tau^{-1}\sigma\tau=\sigma^{4n-1}\rangle\simeq QD_{8n}$
be the quasi-dihedral group of order $16n$.
Let $M$ be the $G$-lattice of rank $4n$ defined in 
{\rm \cite[Definition 7.4]{HKY}}. 
Then $H_{\rm nr}^2(G,M)\simeq \bm{Z}/2\bm{Z}$.
Consequently, $\bm{C}(M)^{G}$ is not retract $\bm{C}$-rational.\\
{\rm (2)}
Let $\widehat G=\langle\sigma^2,\sigma\tau\rangle\simeq Q_{8n}\leq G$ be the
generalized quaternion group of order $8n$.
Let $\widehat M={\rm Res}^G_{\widehat G}(M)$ be the $\widehat G$-lattice
of rank $4n$ defined in {\rm \cite[Definition 7.4]{HKY}}. 
Then $H_{\rm nr}^2(\widehat G,\widehat M)\simeq\bm{Z}/2\bm{Z}\oplus\bm{Z}/2\bm{Z}$.
Consequently,
$\bm{C}(\widehat M)^{\widehat G}$ is not retract $\bm{C}$-rational.
\end{theorem}
\begin{theorem}[{Hoshi, Kang and Yamasaki \cite[Theorem 7.7]{HKY}}]\label{thmCp2p}
Let $p$ be an odd prime and 
$G=\langle \sigma,\tau \mid \sigma^{p^2}=\tau^p=1,
\tau^{-1}\sigma\tau=\sigma^{p+1}\rangle\simeq C_{p^2}\rtimes C_p$.
Let $M$ be the $G$-lattice of rank $p(p-1)$
defined in {\rm \cite[Definition 7.6]{HKY}}. 
Then $H_{\rm nr}^2(G,M)\simeq\bZ/p\bZ$.
Consequently, $\bm{C}(M)^{G}$ is not retract $\bm{C}$-rational.
\end{theorem}


%
%

\end{document}